\documentclass[12 pt]{article}
\usepackage{amsmath,amsfonts,fancybox,multicol,lscape}
\usepackage{color}
\usepackage{graphicx}
\usepackage{tikz}
\usetikzlibrary{arrows,automata,positioning}
\usetikzlibrary{quotes,angles}

\newtheorem{corollary}{Corollary}

\newtheorem{theorem}{Theorem}

\begin{document}

\title{A Direct Proof of the 2nd  Atiyah-Sutcliffe  Conjecture for Convex Quadrilaterals }

\author{Mazen Bou  Khuzam\\
Dept. of Mathematics and Natural Sciences\\
The American University of Iraq at Sulaimani  \\
Email: mazen.boukhuzam@auis.edu.krd}
\date{}
\maketitle

\begin{abstract}
We present a direct proof of the second  conjecture made by M. Atiyah and P. Sutcliffe for the case of convex quadrilaterals. Unlike previous work on this conjecture, our proof does not require any computer aided computations. The new proof relies on a new geometric inequality proved recently by the author.
\end{abstract}

\noindent $\bf{Key words}$: Atiyah determinant ; Atiyah-Sutcliffe conjectures

\begin{figure}[htb!]
 \begin{center}

\end{center}
 \end{figure}

\section{Introduction}

\noindent While studying the spin statistics theorem using classical quantum theory, M. V. Berry and J. M. Robbins \cite{BR} came across the following  purely geometric question in $\Bbb{R}^{3}$ : Does there exist a continuous function from the set of $n$ distinct points in $\Bbb{R}^{3}$ to the flag manifold $U(n)/T^{n}$ which is compatible with the action of the symmetric group? In his attempt to answering their question, M. Atiyah gave a very natural and elegant solution whose validity depends on the nonvanishing of a certain determinant function.  Although  another construction was given by Atiyah in \cite{A1} answering the question of Berry and Robbins in the positive, the nonvanishing of the determinant given in the original  construction remained an unresolved conjecture. Numerical evidence of the validity of the conjecture were provided by Atiyah and  Sutcliffe in their paper \cite{AS} where they added two new conjectures which imply the first one. All three conjectures have remained open since then except for a few verifications on special configurations of points \cite{D} , \cite{EN}, \cite{MP}.

The construction of the determinant begins with $n$ distinct points in $\Bbb{R} ^{3}$. By choosing a prefered axis in $\Bbb{R}^{3}$ we identify it with $\Bbb{R}\times \Bbb{C}$ and denote our $n$ points by $P_{j}=(a_{j},z_{j})$ for $j=1, ... , n$. We consider the vector $V_{jk}=(a_{k}-a_{j}, z_{k}-z_{j})$ to be the vector from $P_{j}$ to $P_{k}$ and use the Hopf map $h:\Bbb{C}^{2} \rightarrow \Bbb{R}^{3}$ given by $h(z,w)=((|z|^{2}-|w|^{2})/2 , z\overline{w} )$ to lift $V_{jk}$ from $\Bbb{R}^{3}$ to $\Bbb{C}^{2}$.  If we denote also $a_{k}-a_{j}$ by $a_{jk}$,  $z_{k}-z_{j}$ by $z_{jk}$ and set $d_{jk}=a_{jk}+\sqrt{a_{jk}^{2}+|z_{jk}|^{2}}$  we see that the lift of $V_{jk}$ is given by 

\[  h^{-1}  (V_{jk}) =  e^{i  \theta } \frac{1}{\sqrt{d_{jk}}} (d_{jk} , \overline{z}_{jk} )
\hspace{1.5 cm} ; \, \, \theta \in \Bbb{R}, \, i=\sqrt{-1} \] 


Since $h^{-1}  (V_{jk})$ does not uniquely depend on $V_{jk}$ we  follow the normalization of Atiyah by choosing $(-\overline{w}, \overline{z})$ as a lift (under $h$) for $V_{kj}$ whenever $(z,w)$ is a lift for  $V_{jk}$ (see \cite{A2}). We let $C_{k}$ be the symmetric tensor product of $ h^{-1}(V_{kj})$  for all $j\neq k$. The determinant of the $n\times n$ matrix whose $k$th column is $C_{k}$ will be called the Atiyah  determinant and will be denoted by $At(P_{1},...,P_{n})$ or just $At$ if the points are already specified. The determinant function $At$  is invariant under rotations and translations and gets conjugated under plane reflections. The first  conjecture (C1) of Atiyah is $At \neq 0$ while the second Atiyah-Sutcliffe conjecture (C2) is $\displaystyle{ At \geq \prod _{j<k} (2 || \overrightarrow{P_{j}P_{k}}||)}$, and the third one (C3) is

\[ At^{n-2}(P_{1}, ... , P_{n}) \geq \prod _{j=1}^{n} At(P_{1}, ... , P_{j-1}, P_{j+1}, P_{n}) \]

\noindent for all  distinct points $P_{1}$, ..., $P_{n}$ of $\Bbb{R} ^{3}$. It is clear from the statement of these conjectures that $(C3)\implies (C2) \implies (C1)$.

Most of the methods applied in previous work on this determinant considered special cases of configurations of points (see for example \cite{D}, \cite{MP}) or used some sort of computer aid in their computations. The main work on the general four-point case was done by M. Eastwood and P. Norbury in \cite{EN} where they  proved the first conjecture of Atiyah to be true for $n=4$ points in space.  Others followed their formula using Maple (see \cite{M1}) and proved the three conjectures with computer aid. Our paper provides the first direct proof of the second Atiyah-Sutcliffe conjecture for four  points forming a convex quadrilateral. The proof we present here does not rely on any computer calculations, but rather on a geometric inequality proved by the author in a previous paper \cite{M}. We hope that this new approach will pave the way for generalizations that may solve the full conjecture.


\section{The Planar Case}

In this section, we explore the Atiyah determinant for four planar points and present it in standard form. First note that $At$ is a real number when all points lie in a plane. This is because a reflection in their plane  leaves the points fixed since $\overline{At}=At$.  This was already mentioned in \cite{A2} and will also be seen here concretely from the expansion of its standard form. Moreover, after applying a solid motion, we may assume  our points are in $\{ 0 \} \times \Bbb{C}$ so that all $a_{jk} $ values  
become zeros and our points can be written as  $(0,z_{1})$, ..., 
$(0,z_{4})$. Based on that, the construction of the Atiyah determinant is done as follows:

When the first point is considered as an observer of the three other points we obtain 
$(0,z_{12})$, $(0,z_{13})$, $(0,z_{14})$ and these lift under the Hopf map $h$  to 
\[ \frac{1}{\sqrt{r_{12}}} (r_{12}, \overline{z} _{12}), \, \, \, \, \, \, \,
\frac{1}{\sqrt{r_{13}}} (r_{13}, \overline{z} _{13}), \, \, \, \, \, \, \, 
\frac{1}{\sqrt{r_{14}}} (r_{14}, \overline{z} _{14})  \]
 
\noindent where $r_{jk}$ denotes $|z_{jk}|$.

%

\noindent We do the same thing when $(0,z_{2})$ is a vision point and get the vectors
 
\noindent $(0,z_{21})$, $(0,z_{23})$, $(0,z_{24})$ whose lifts are 
\[ \frac{1}{\sqrt{r_{21}}} (z _{21} , r_{21} ), \, \, \, \, \, \, \,
\frac{1}{\sqrt{r_{23}}} (r_{23}, \overline{z} _{23}), \, \, \, \, \, \, \, 
\frac{1}{\sqrt{r_{24}}} (r_{24}, \overline{z} _{24}).  \]

\noindent Similarly, the lifts corresponding for the vision points $(0,z_{3})$ and $(0,z_{4})$ are

\[ \frac{1}{\sqrt{r_{31}}} (z _{31} , r_{31} ), \, \, \, \, \, \, \,
\frac{1}{\sqrt{r_{32}}} (z_{32}, r _{32}), \, \, \, \, \, \, \, 
\frac{1}{\sqrt{r_{34}}} (r_{34}, \overline{z} _{34}).  \]

\noindent and

\[ \frac{1}{\sqrt{r_{41}}} (z _{41} , r_{41} ), \, \, \, \, \, \, \,
\frac{1}{\sqrt{r_{42}}} (z_{42}, r _{42}), \, \, \, \, \, \, \, 
\frac{1}{\sqrt{r_{43}}} (z_{43}, r _{43}).  \]

\noindent respectively.

\noindent Using $u_{jk}=\frac{z_{jk}}{r_{jk}}$ to be the direction of $z_{jk}$ we see that the lifts can be written as: 

\[ \sqrt{r_{12}} (1,\overline{u} _{12}), \, \, \, \, \, \, 
\sqrt{r_{13}} (1,\overline{u} _{13}), \, \, \, \, \, \, 
\sqrt{r_{14}} (1,\overline{u} _{14}) \]

\[ \sqrt{r_{21}} u_{21} (1,\overline{u} _{21}), \, \, \, \, \, \, 
\sqrt{r_{23}} (1,\overline{u} _{23}), \, \, \, \, \, \, 
\sqrt{r_{24}} (1,\overline{u} _{24}) \]

\[ \sqrt{r_{31}} u_{31} (1,\overline{u} _{31}), \, \, \, \, \, \, 
\sqrt{r_{32}} u_{32} (1,\overline{u} _{32}), \, \, \, \, \, \, 
\sqrt{r_{34}} (1,\overline{u} _{34}) \]

\[ \sqrt{r_{41}} u_{41} (1,\overline{u} _{41}), \, \, \, \, \, \, 
\sqrt{r_{42}} u_{42} (1,\overline{u} _{42}), \, \, \, \, \, \, 
\sqrt{r_{43}} u_{43} (1,\overline{u} _{43}) \]

\noindent The symmetric tensor product of the vectors in each of the lines above gives us: 

\[ \sqrt{r_{12} r_{13} r_{14}} (1, 
\overline{u}_{12}+\overline{u}_{13}+\overline{u}_{14}, \, \, \, \, \, \, 
\overline{u}_{12}\overline{u}_{13}+\overline{u}_{12}\overline{u}_{14}+
\overline{u}_{13}\overline{u}_{14}, \, \, \, \, \, \, 
\overline{u}_{12}\overline{u}_{13}\overline{u}_{14} ) \] 

\[ \sqrt{r_{21} r_{23} r_{24}} u_{21} (1, 
\overline{u}_{21}+\overline{u}_{23}+\overline{u}_{24}, \, \, \, \, \, \,  
\overline{u}_{21}\overline{u}_{23}+\overline{u}_{21}\overline{u}_{24}+
\overline{u}_{23}\overline{u}_{24}, \, \, \, \, \, \, 
\overline{u}_{21}\overline{u}_{23}\overline{u}_{24} ) \]

\[ \sqrt{r_{31} r_{32} r_{34}} u_{31} u_{32} (1, 
\overline{u}_{31}+\overline{u}_{32}+\overline{u}_{34}, \, \, \, \, \, \,  
\overline{u}_{31}\overline{u}_{32}+\overline{u}_{31}\overline{u}_{34}+
\overline{u}_{32}\overline{u}_{34}, \, \, \, \, \, \, 
\overline{u}_{31}\overline{u}_{32}\overline{u}_{34} ) \]

\[ \sqrt{r_{41} r_{42} r_{43}} u_{41} u_{42} u_{43} (1, 
\overline{u}_{41}+\overline{u}_{42}+\overline{u}_{43}, \, \, \, \, \, \,  
\overline{u}_{41}\overline{u}_{42}+\overline{u}_{41}\overline{u}_{43}+
\overline{u}_{42}\overline{u}_{43}, \, \, \, \, \, \, 
\overline{u}_{41}\overline{u}_{42}\overline{u}_{43} ). \]

\noindent Consequently, if we define the angular part of $At$ to be $At_{ang}$ given by the product   $u_{12}u_{13}u_{14}u_{23}u_{24}u_{34}$ times the determinant

\begin{center} $   \left| \begin{array}{cccc}
                    1 & 1 & 1 & 1   \\
                                    \\

 \overline{u} _{12}+ \overline{u} _{13} 

 & \overline{u} _{21}+ \overline{u} _{23} & \overline{u} _{31}+ \overline{u} _{32}  & \overline{u} _{41}+ \overline{u} _{42}    \\
+ \overline{u} _{14} &  + \overline{u} _{24} & + \overline{u} _{34} & + \overline{u} _{43} \\

                                                 \\
                    \overline{u} _{12}\overline{u} _{13}+ \overline{u} _{12}\overline{u} _{14} & \overline{u} _{21}\overline{u} _{23}+ \overline{u} _{21}\overline{u} _{24} & \overline{u} _{31}\overline{u} _{32}+ \overline{u} _{31}\overline{u} _{34} & \overline{u} _{41}\overline{u} _{42}+ \overline{u} _{41}\overline{u} _{43}   \\

+ \overline{u} _{13}\overline{u} _{14} & + \overline{u} _{23}\overline{u} _{24} & + \overline{u} _{32}\overline{u} _{34} & + \overline{u} _{42}\overline{u} _{43} \\
 
\\
                  \overline{u} _{12}\overline{u} _{13}\overline{u} _{14}  & \overline{u} _{21}\overline{u} _{23}\overline{u} _{24} & \overline{u} _{31}\overline{u} _{32}\overline{u} _{34} & \overline{u} _{41}\overline{u} _{42}\overline{u} _{43}   \\
                  
          \end{array} \right|$ \end{center}

\noindent then  the Atiyah determinant will be $At=P At_{ang}$  where $P= r_{12}r_{13}r_{14}r_{23}r_{24}r_{34}$ is the product of the lengths of all sides. The decomposition $At=PAt_{ang}$ will be called the standard form of the Atiyah determinant for four planar points. It is also natural to call $P$ the scalar part of $At$ since the rest of the product $At_{ang}$ is the angular part of $At$.

\section{The Convex Case}

\begin{figure}[htb!]

\centering

\begin{tikzpicture}
[scale=2, every node/.style={anchor=center, scale=1.5}]
  \draw
    (-2,0) coordinate (a) node[right] {}
    (4,0) coordinate (b) node[left] {}
    (1,-1.32) coordinate (c) node[above right] {}
    (0.7,-1.32) coordinate (z4) node[above right] {{\footnotesize $z_1$}}
    (-0.5,-2) coordinate (d) node[above right] {}
    (-1.1,-2.2) coordinate (z1) node[above right] {{\footnotesize $z_2$}}
    (2.5,-2) coordinate (e) node[above right] {}
     (0.5,-2.4) coordinate (z1) node[above right] {{\footnotesize $W'$}}
    (2.5,-2) coordinate (e) node[above right] {}
    (0.95,-2.4) coordinate (z1) node[above right] {{\footnotesize $W$}}
    (2.5,-2) coordinate (e) node[above right] {}
  (-0.45,-2.35) coordinate (z1) node[above right] {{\footnotesize $\alpha _{2}$}}
 (-0.2,-2.1) coordinate (z1) node[above right] {{\footnotesize $\beta _{2}$}}
    (2.5,-2) coordinate (e) node[above right] {}
  (1.8,-2.1) coordinate (z1) node[above right] {{\footnotesize $\alpha _{4}$}}
    (2.5,-2) coordinate (e) node[above right] {}
  (1.9,-2.45) coordinate (z1) node[above right] {{\footnotesize $\beta _{4}$}}
    (2.5,-2) coordinate (e) node[above right] {}
(0.95,-3.7) coordinate (z1) node[above right] {{\footnotesize $\alpha _{3}$}}
    (2.1,-2) coordinate (e) node[above right] {}
  (0.6,-3.7) coordinate (z1) node[above right] {{\footnotesize $\beta _{3}$}}
    (2.5,-2) coordinate (e) node[above right] {}
    (2.5,-2) coordinate (e) node[above right] {}
    (2.5,-2.2) coordinate (z3) node[above right] {{\footnotesize $z_4$}}
    (1,-4) coordinate (f) node[above right] {}
    (0.7,-4.5) coordinate (z2) node[above right] {{\footnotesize $z_3$}};  
 \draw[dashed]  (-2,0) -- (1,-1.32);
 \draw   (1,-1.32) -- (2.5,-2);
   \draw[dashed] (-0.5,-2) -- (-2,0);
\draw (-0.5, -2) -- (1,-4);
\draw[dashed] (4,0) -- (2.5,-2);  
\draw (2.5,-2) -- (1,-4);
\draw[dashed]    (4,0) -- (1,-1.32);
\draw   (1,-1.32) -- (-0.5,-2)
    (-0.5,-2) -- (2.5,-2)
    (1,-1.32) -- ((1,-4)
    pic["$X$", angle eccentricity=1.5, angle radius=0.5cm]
    {angle=f--a--e}
    pic["$Y$", angle eccentricity=1.5, angle radius=0.5cm]
    {angle=d--b--f}
    pic["{\footnotesize $\alpha_1$}", angle eccentricity=1.7, angle radius=0.25cm]
    {angle=d--c--f}
    pic["{\footnotesize $\beta_1$}", angle eccentricity=1.7, angle radius=0.25cm]
    {angle=f--c--e};
\end{tikzpicture}
 
\caption{Convex Quadrilateral}
\label{fig1}
\end{figure}
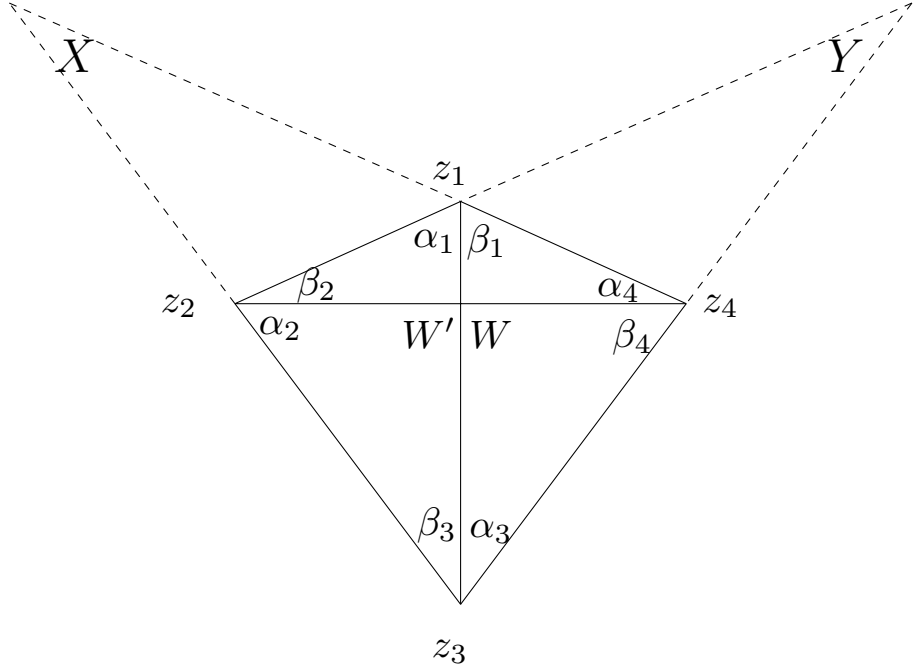

\noindent Let us consider the terms of the Atiyah determinant corresponding to the permutation $(1,2,3,4)$ in the expansion of $At$. This corresponds to the product of the main diagonal entries: 

\[ ( \overline{u} _{21}+ \overline{u} _{23} + \overline{u} _{24} )( \overline{u} _{31}\overline{u} _{32}+ \overline{u} _{31}\overline{u} _{34} + \overline{u} _{32}\overline{u} _{34}) \overline{u} _{41}\overline{u} _{42}\overline{u} _{43} \]

\vspace{2 mm}

\noindent multiplied with $r_{12}r_{13}r_{14}r_{23}r_{24}r_{34}u_{12}u_{13}u_{14}u_{23}u_{24}u_{34}$. We decompose this expression as an angular part $ ( 1 + u _{21} \overline{u} _{23} +  u _{21} \overline{u} _{24} )( 1 +  u _{31}\overline{u} _{34} + u_ {32} \overline{u} _{34})$ multiplied with the scalar part $ P = r_{12}r_{13}r_{14}r_{23}r_{24}r_{34}$. The angular part, expressed with respect to the angles appearing in figure 1, can be written as 

\[ (1+e^{i\gamma _{2}} + e^{i \beta _{2}})(1+ e^{i \alpha _{3} } + e^{i \gamma _{3}}) \]

 This presentation of the terms of the Atiyah determinant can be done in a similar way for any permutation of the set $\{ 1, 2, 3, 4 \}$. For example, the terms corresponding to the permutation $(1,4,3,2)$ can also be written as the scalar part  $P=  r_{12}r_{13}r_{14}r_{23}r_{24}r_{34}$ multiplied with the angular part $(1+e^{-i\gamma _{4}} + e^{-i\alpha _{4}})(1+e^{-i \beta _{3} } + e^{-i \gamma _{3}})$. Taking all terms of the Atiyah determinant into account requires considering all 24 permutations. Let us first  start with the terms of the identity permutation: Expanding the angular part corresponding for the identity permutation yields a linear part  given by  $\displaystyle{ 1 + e^{i\gamma _{2}} +  e^{i\beta _{2}} + e^{i\gamma _{3}} + e^{i\alpha _{3}} }$ and a quadratic part given by $e^{i(\gamma _{2}+\gamma _{3})} +  e^{i(\beta _{2}+ \gamma _{3})} +  e^{i(\gamma _{2}+\alpha _{3})} +  e^{i(\beta _{2}+ \alpha _{3})}$. We need to  collect all similar terms from all permutations to get a general formula for the Atiyah determinant:

\underline{The Linear Terms}: These terms can be summarized as follows: We have 24 ones, one from each permutation. We also have all angles of the quadrilaterial appearing with an equal number of times as $e^{i \theta}$, where $\theta$ is any one of the $\alpha _{j}$'s, $\beta _{j}$'s or $\gamma _{j}$'s, and  $j=1,2,3,4$. It is easy to see that conjugate terms appear from permutations of the form $(a,b,c,d)$ and their reverse order $(d,c,b,a)$, which confirms that the determinant is a real number. This allows us to see that each linear term $e^{i\theta}$ appears 4 times and each $e^{-i\theta}$ appears 4 times. Accordingly, the linear terms in the expansion of the Atiyah determinant are 

\[ 24 + 8 \sum  \cos \theta _{k} \]

 where $\theta _{k}$ runs over  the  12 angles $\alpha _{j}$, $\beta _{j}$, and $\gamma _{j}$, where $j=1,2,3,4$. 

\vspace{2 mm}

\underline{The Quadratic Terms}:  Listing all quadratic terms based on their corresponding permutations gives us: 


\vspace{2 mm}

\noindent (1,2,3,4): $\displaystyle{e^{i(\gamma _{2}+\gamma _{3})} +  e^{i(\beta _{2}+ \gamma _{3})} +  e^{i(\gamma _{2}+\alpha _{3})} +  e^{i(\beta _{2}+ \alpha _{3})}}$

\vspace{1 mm}

\noindent (1,2,4,3): $\displaystyle{e^{i(\beta _{2}-\beta _{4})} +  e^{i(\gamma _{2}- \beta _{4})} +  e^{i(\beta _{2}-\gamma _{4})} +  e^{i(\gamma _{2}- \gamma _{4})}}$

\vspace{1 mm}

\noindent (1,3,2,4): $\displaystyle{e^{i(-\beta _{3}-\alpha _{2})} +  e^{i(\alpha _{3}- \alpha _{2})} +  e^{i(-\beta _{3}+\beta _{2})} +  e^{i(\alpha _{3}+ \beta _{2})}}$ 

\vspace{1 mm} 

\noindent (1,3,4,2): $\displaystyle{e^{i(\alpha _{3}+\beta _{4})} +  e^{i(-\beta _{3}+ \beta _{4})} +  e^{i(\alpha _{3}-\alpha _{4})} +  e^{i(-\beta _{3}- \alpha _{4})}}$ 

\vspace{1 mm}

\noindent (1,4,2,3): $\displaystyle{e^{i(-\alpha _{4}+\alpha _{2})} +  e^{i(-\gamma _{4}+ \alpha _{2})} +  e^{i(-\alpha _{4}+\gamma _{2})} +  e^{i(-\gamma _{4}+ \gamma _{2})}}$

\noindent (1,4,3,2): $\displaystyle{e^{i(-\gamma _{4}-\gamma _{3})} +  e^{i(-\alpha _{4}- \gamma _{3})} +  e^{i(-\gamma _{4}-\beta _{3})} +  e^{i(-\alpha _{4}- \beta _{3})}}$ 

\vspace{1 mm}

\noindent (2,1,3,4): $\displaystyle{e^{i(-\alpha _{1}+\alpha _{3})} +  e^{i(-\gamma _{1}+ \alpha _{3})} +  e^{i(-\alpha _{1}+\gamma _{3})} +  e^{i(-\gamma _{1}+ \gamma _{3})}}$ 

\vspace{1 mm}

\noindent (2,1,4,3): $\displaystyle{e^{i(-\gamma _{1}-\gamma _{4})} +  e^{i(-\alpha _{1}- \gamma _{4})} +  e^{i(-\gamma _{1}-\beta _{4})} +  e^{i(-\alpha _{1}- \beta _{4})}}$ 

\vspace{1 mm}

\noindent (3,1,2,4): $\displaystyle{e^{i(\alpha _{1}+\beta _{2})} +  e^{i(-\beta _{1}+ \beta _{2})} +  e^{i(\alpha _{1}-\alpha _{2})} +  e^{i(-\beta _{1}- \alpha _{2})}}$ 

\vspace{1 mm}

\noindent (3,1,4,2): $\displaystyle{e^{i(-\beta _{1}-\alpha _{4})} +  e^{i(\alpha _{1}- \alpha _{4})} +  e^{i(-\beta _{1}+\beta _{4})} +  e^{i(\alpha _{1}+ \beta _{4})}}$ 

\vspace{1 mm}

\noindent (4,1,2,3): $\displaystyle{e^{i(\gamma _{1}+\gamma _{2})} +  e^{i(\beta _{1}+ \gamma _{2})} +  e^{i(\gamma _{1}+\alpha _{2})} +  e^{i(\beta _{1}+ \alpha _{2})}}$ 

\vspace{1 mm}

\noindent (4,1,3,2): $\displaystyle{e^{i(\beta _{1}-\beta _{3})} +  e^{i(\gamma _{1}- \beta _{3})} +  e^{i(\beta _{1}-\gamma _{3})} +  e^{i(\gamma _{1}- \gamma _{3})}}$

\vspace{3 mm}

 The other 12 permutations are nothing but the conjugates of the ones appearing above. Accordingly, adding all terms would result in a sum of twice the cosines of the angles appearing in this list. We note that all the terms appearing in the first column will cancel out. For example, $\cos (\gamma _{2} + \gamma _{3}) + \cos (\alpha _{3} - \alpha _{1})=0$. This can be seen to apply to all terms of the first column. On the other hand, the terms of the fourth column appear twice in the list (including the conjugates). 

Summarizing all terms can be done as follows: the multiplicity-one terms can be put in 6 families:   $2 \cos (\alpha _{j}  -\alpha _{j+1})$, $2 \cos (\beta _{j} -\beta _{j+1})$, $2 \cos (\gamma _{j} +\alpha _{j+1})$, $2 \cos (\gamma _{j+1} +\beta _{j})$, $2 \cos (\gamma _{j} -\alpha _{j+2})$, $2 \cos (\gamma _{j} -\beta _{j+2})$, where  $j=1,2,3,4$, modulo 4, and the multiplicity-two terms can be listed as follows: $4 \cos (\alpha _{1} + \beta _{4})$, $4 \cos (\alpha _{2} + \beta _{1})$,  $4 \cos (\alpha _{3} + \beta _{2})$, $4 \cos (\alpha _{4} + \beta _{3})$, $4 \cos (\gamma _{1} - \gamma _{3})$, and $4 \cos (\gamma _{2} - \gamma _{4})$. Contemplating this list, we realize that these terms can be grouped based on vertices facing opposite triangles in the following manner:

\vspace{1 cm}

\noindent $ \cos (\alpha _{j}  -\alpha _{j+1})$:   where $\alpha _{j}$ is at vertex $z_{j}$ and $\alpha _{j+1}$ is in the opposite triangle

\vspace{2 mm}

\noindent $ \cos (\beta _{j} -\beta _{j-1})$:  where $\beta _{j}$ is at vertex $z_{j}$ and $\beta _{j-1}$ is in the opposite triangle

\vspace{2 mm} 

\noindent $ \cos (\gamma _{j} +\alpha _{j+1})$:  where $\gamma _{j}$ is at vertex $z_{j}$ and $\alpha _{j+1}$ is in the opposite triangle

\vspace{2 mm} 

\noindent $\cos (\gamma _{j} +\beta _{j-1})$:  where $\gamma _{j}$ is at vertex $z_{j}$ and $\beta _{j-1}$ is in the opposite triangle

\vspace{2 mm} 

\noindent $ \cos (\alpha _{j} -\gamma _{j+2})$:  where $\alpha _{j}$ is at vertex $z_{j}$ and $\gamma _{j+2}$ is in the opposite triangle

\vspace{2 mm} 

\noindent $ \cos (\beta _{j} -\gamma _{j+2})$:  where $\beta _{j}$ is at vertex $z_{j}$ and $\gamma _{j+2}$ is in the opposite triangle

\vspace{2 mm}

\noindent in addition to the multiplicity-two terms:

\vspace{3 mm}

\noindent $ \cos (\alpha _{1} + \beta _{4})$: where $\alpha _{1}$ is at vertex $z_{1}$ and $\beta _{4}$ is in the opposite triangle

\hspace{1.8 cm} or   $\beta _{4}$ is at vertex $z_{4}$ and $\alpha _{1}$ is in the opposite triangle

\noindent $\cos (\alpha _{2} + \beta _{1})$: where $\alpha _{2}$ is at vertex $z_{2}$ and $\beta _{1}$ is in the opposite triangle

\hspace{1.8 cm} or   $\beta _{1}$ is at vertex $z_{1}$ and $\alpha _{2}$ is in the opposite triangle

\noindent $ \cos (\alpha _{3} + \beta _{2})$: where $\alpha _{3}$ is at vertex $z_{3}$ and $\beta _{2}$ is in the opposite triangle

\hspace{1.8 cm} or   $\beta _{2}$ is at vertex $z_{2}$ and $\alpha _{3}$ is in the opposite triangle

\noindent $ \cos (\alpha _{4} + \beta _{3})$: where $\alpha _{4}$ is at vertex $z_{4}$ and $\beta _{3}$ is in the opposite triangle

\hspace{1.8 cm} or   $\beta _{3}$ is at vertex $z_{3}$ and $\alpha _{4}$ is in the opposite triangle

\noindent $ \cos (\gamma _{1} - \gamma _{3})$: where $\gamma _{1}$ is at vertex $z_{1}$ and $\gamma _{3}$ is in the opposite triangle

\hspace{1.8 cm} or   $\gamma _{3}$ is at vertex $z_{3}$ and $\gamma _{1}$ is in the opposite triangle

\noindent $ \cos (\gamma _{2} - \gamma _{4})$: where $\gamma _{2}$ is at vertex $z_{2}$ and $\gamma _{4}$ is in the opposite triangle

\hspace{1.8 cm} or   $\gamma _{4}$ is at vertex $z_{4}$ and $\gamma _{2}$ is in the opposite triangle

\vspace{3 mm}

 For example, let us consider the vertex $z_{1}$ facing the triangle $z_{2}z_{3}z_{4}$. The angle $\alpha _{1}$ at the vertex $z_{1}$ is coupled with each of the angles of the triangle  $z_{2}z_{3}z_{4}$ in the terms: $\cos (\alpha _{1} - \alpha _{2})$, $\cos (\alpha _{1} - \gamma _{3})$, and $\cos (\alpha _{1} + \beta _{4})$. Also,  the angle $\beta _{1}$ at the vertex $z_{1}$  is coupled with each of the angles of the triangle  $z_{2}z_{3}z_{4}$ in the terms: $\cos (\beta _{1} + \alpha _{2})$, $\cos (\beta _{1} - \gamma _{3})$, and $\cos (\beta _{1} - \beta _{4})$. Finally, the angle $\gamma _{1}$  is coupled with each of the angles of the triangle  $z_{2}z_{3}z_{4}$ in the terms: $\cos (\gamma _{1} + \alpha _{2})$, $\cos (\gamma _{1} - \gamma _{3})$, and $\cos (\gamma _{1} + \beta _{4})$.

Collecting the product of cosines $\cos \alpha \cos \beta$ from $\cos (\alpha \pm \beta ) = \cos \alpha \cos \beta \mp \sin \alpha \sin \beta $ for each of the quadratic terms we obtain the sum of products

\[  2  \sum _{j=1}^{4} \left( \cos \alpha _{j} +  \cos \beta _{j} + \cos \gamma _{j} \right) \left( \cos \alpha _{j+1} + \cos \gamma _{j+2} + \cos \beta _{j+3} \right) \]

\vspace{2 mm}

\noindent where all indices are taken modulo four. Here, we see that $\displaystyle{ \cos \alpha _{j} +  \cos \beta _{j}}$

\noindent $\displaystyle{ + \cos \gamma _{j} }$ is the sum of cosines of all angles at the vertex $z_{j}$  
\noindent whereas $ \cos \alpha _{j+1} + \cos \gamma _{j+2} + \cos \beta _{j+3}  $ is the sum of cosines of all angles of the opposite triangle.

To accommodate for the  remaining product of sines $\mp \sin \alpha \sin \beta$ from $\cos (\alpha \pm \beta ) = \cos \alpha \cos \beta \mp \sin \alpha \sin \beta $, let us call the sum of all products of two sines (together with their signs) to be the expression  $E$. Based on that we have proved the following theorem:

\begin{theorem}: The Atiyah determinant can be written as $r_{12}r_{13}r_{14}r_{23}r_{24}r_{34}$ times the angular part

\[ 24 + 8 \sum _{j=1}^{4}  \left( \cos \alpha _{j} +  \cos \beta _{j} +  \cos \gamma _{j} \right) \]

 \[ + \, 2  \sum _{j=1}^{4} \left( \cos \alpha _{j} +  \cos \beta _{j} + \cos \gamma _{j} \right) \left( \cos \alpha _{j+1} + \cos \gamma _{j+2} + \cos \beta _{j+3} \right)  \]

\[ + E \]

\end{theorem}

\vspace{4 mm}

\noindent Let us call the first sum $\displaystyle{ S_{1} = \sum _{j=1}^{4}  \left( \cos \alpha _{j} +  \cos \beta _{j} +  \cos \gamma _{j} \right)}$ and the second sum $\displaystyle{S_{2}=  \sum _{j=1}^{4} \left( \cos \alpha _{j} +  \cos \beta _{j} + \cos \gamma _{j} \right) \left( \cos \alpha _{j+1} + \cos \gamma _{j+2} + \cos \beta _{j+3} \right)  }$

\vspace{3 mm}

\begin{corollary}:  \[ S_{1} = \sum _{j=1}^{4}  \left( \cos \alpha _{j} +  \cos \beta _{j} +  \cos \gamma _{j} \right) \geq 4 \]

\end{corollary}

\vspace{2 mm}

\noindent Proof: This proof relies on the identity 

\[ \cos \alpha + \cos \beta + \cos \gamma = 1 + 4 \sin \frac{\alpha }{2} \sin \frac{\beta }{2} \sin \frac{\gamma }{2} \hspace{5 mm} ; \hspace{5 mm}  \alpha + \beta + \gamma = \pi \]

\hfill (*)

Accordingly, we regroup the terms in our summation $S_{1}$ in such a way  that angles of the same triangle are grouped together to obtain

\[ S_{1} =  \sum _{j=1}^{4}  \left( \cos \alpha _{j+1} + \cos \gamma _{j+2} + \cos \beta _{j+3} \right) \] 

\[ = \sum _{j=1}^{4}  \left( 1+ 4 \sin \frac{ \alpha _{j+1}}{2} \sin \frac{\gamma _{j+2}}{2}  \sin \frac{ \beta _{j+3}}{2} \right) \] 

\[ = 4 + 4 \sum _{j=1}^{4}   \sin \frac{ \alpha _{j+1}}{2} \sin \frac{\gamma _{j+2}}{2}  \sin \frac{ \beta _{j+3}}{2}   \] 

\noindent which is obviously greater than or equal to 4 since all other terms are non-negative.

\vspace{4 mm}

\begin{corollary}:  \[ 8 S_{1} + 2 S_{2} \geq 40 \]

\end{corollary}

\noindent Proof: First, we collect from corollary 1 the simplified expression of $S_{1}$ and write it as

\[ 8 S_{1} = 32 + 32 \sum _{j=1}^{4}   \sin \frac{ \alpha _{j+1}}{2} \sin \frac{\gamma _{j+2}}{2}  \sin \frac{ \beta _{j+3}}{2}   \] 

\vspace{2 mm}

We now turn our attention to $S_{2}$. Here we use the identity (*) again on the second factors of the sum $S_{2}$ to write it as

\[ S_{2}=  \sum _{j=1}^{4} \left( \cos \alpha _{j} +  \cos \beta _{j} + \cos \gamma _{j} \right) (  1 + 4  \sin \frac{ \alpha _{j+1}}{2} \sin \frac{\gamma _{j+2}}{2}  \sin \frac{ \beta _{j+3}}{2}   )  \]

\vspace{2 mm}

Consequently, 

\[ 2 S_{2}=  2 \sum _{j=1}^{4} \left( \cos \alpha _{j} +  \cos \beta _{j} + \cos \gamma _{j} \right) \]

\[ + 8 \sum _{j=1}^{4}    \left( \cos \alpha _{j} +  \cos \beta _{j} + \cos \gamma _{j} \right)   \sin \frac{ \alpha _{j+1}}{2} \sin \frac{\gamma _{j+2}}{2}  \sin \frac{ \beta _{j+3}}{2}   \]

The first line of $2S_{2}$ can also be simpified since angles can be regrouped  into four groups, each corresponding to one of the four triangles of the quadrilateral, so we write as 

\[  2 \sum _{j=1}^{4} \left( \cos \alpha _{j}  + \cos \gamma _{j+1} + \cos \beta _{j+2} \right) = \]

\[ = 8 + 8 \sum _{j=1}^{4}   \sin \frac{ \alpha _{j+1}}{2} \sin \frac{\gamma _{j+2}}{2}  \sin \frac{ \beta _{j+3}}{2}   \] 

\vspace{2 mm}

Accordingly, $8S_{1}+2S_{2}$ is equal to

\[ = 40 + 40 \sum _{j=1}^{4}   \sin \frac{ \alpha _{j+1}}{2} \sin \frac{\gamma _{j+2}}{2}  \sin \frac{ \beta _{j+3}}{2}   \] 

\[ + 8 \sum _{j=1}^{4}    \left( \cos \alpha _{j} +  \cos \beta _{j} + \cos \gamma _{j} \right)   \sin \frac{ \alpha _{j+1}}{2} \sin \frac{\gamma _{j+2}}{2}  \sin \frac{ \beta _{j+3}}{2}   \]

\vspace{2 mm}

which finally gives us that 

\[ 8S_{1}+2S_{2} = 40 + 8 \sum _{j=1}^{4}    \left( 5+ \cos \alpha _{j} +  \cos \beta _{j} + \cos \gamma _{j} \right)   \sin \frac{ \alpha _{j+1}}{2} \sin \frac{\gamma _{j+2}}{2}  \sin \frac{ \beta _{j+3}}{2}   \]

Since  $5+ \cos \alpha _{j} +  \cos \beta _{j} + \cos \gamma _{j} \geq 0$, we can easily see that $8S_{1}+2S_{2}\geq 40$ which proves the second corollary.

\vspace{2 mm} 

Based on these two corollaries, we can prove the second Atiyah-Sutcliffe conjecture for a convex quadrilateral if we can prove that $E\geq 0$. This will finish the proof. To do that, we will prove that $E$ is nothing but the expression $4(E_{12} + E_{23} + E_{34} + E_{14} - E_{13} - E_{24}) $ which was proved to be non-negative in \cite{M}. Before we start comparing our work with \cite{M}, note that the side-lengths in our quadrilateral were denoted in \cite{M} as follows: $r_{23}=a$, $r_{13}=b$, $r_{12}=c$, $r_{14}=d$, $r_{24}=e$, and $r_{34}=f$. Let us proceed by  regrouping the terms in the following way:

\vspace{2 mm}

Consider each side of the quadrilaterial separately. For example, let us consider the products of sines corresponding to the side $z_{3}z_{4}$. This side belong to two triangles, namely $z_{2}z_{3}z_{4}$ and $z_{1}z_{3}z_{4}$. When $z_{1}$ is facing triangle $z_{2}z_{3}z_{4}$, the terms corresponding for $z_{3}z_{4}$ are all the terms corresponding for the angle $\beta _{1}$ at $z_{1}$ combined with the angles of the triangle $z_{2}z_{3}z_{4}$. These are: $\cos (\beta _{1}-\beta _{4})$,  $\cos (\beta _{1}-\gamma _{3})$, and $\cos (\beta _{1}+\alpha _{2})$. By doing a similar choice of the terms when $z_{2}$ is facing triangle  $z_{1}z_{3}z_{4}$, we see that the terms corresponding for the side  $z_{3}z_{4}$ are: $\cos (\alpha _{2}-\alpha _{3})$,  $\cos (\alpha _{2}-\gamma _{4})$, and $\cos (\alpha _{2}+\beta _{1})$. When we collect the product of two sines from these terms as was formed in $E$ we obtain

\[ \sin\beta _{1} \sin\beta _{4} + \sin \beta _{1} \sin \gamma _{3}  + \sin \alpha _{2} \sin \alpha _{3} + \sin \alpha _{2} \sin \gamma _{4} - 2  \sin \alpha _{2} \sin \beta _{1} \]

\vspace{2 mm}

Multiplying this with the scalar part $abcdef$, we note that we can write this expression in terms of the areas of the triangles as follows: 

\vspace{2 mm}

$\displaystyle{(abcdef)  \sin\beta _{1} \sin\beta _{4} = 4 ac (\frac{1}{2} bd \sin\beta _{1} ) (\frac{1}{2} ef  \sin\beta _{4}) =4ac A_{134}A_{234}}$

\vspace{1 mm}

$\displaystyle{(abcdef)  \sin \beta _{1} \sin \gamma _{3} = 4 ce (\frac{}{2} bd \sin\beta _{1} ) (\frac{1}{2}  af \sin \gamma _{3}) =4ce A_{134}A_{234}}$

\vspace{1 mm}

$\displaystyle{(abcdef)  \sin \alpha _{2} \sin \alpha _{3} = 4 cd (\frac{}{2} ae \sin \alpha _{2}  ) (\frac{1}{2}  bf \sin \alpha _{3}) =4cd A_{134}A_{234}}$

\vspace{1 mm}

$\displaystyle{(abcdef) \sin \alpha _{2} \sin \gamma _{4} = 4 bc (\frac{}{2} ae \sin \alpha _{2}  ) (\frac{1}{2} df  \sin \gamma _{4}) =4bc A_{134}A_{234}}$

\vspace{1 mm}

$\displaystyle{  (abcdef)   \sin \alpha _{2} \sin \beta _{1} = 4 cf (\frac{}{2} ae \sin \alpha _{2}  ) (\frac{1}{2} bd  \sin \beta _{1}) =4cf A_{134}A_{234}}$

\vspace{2 mm}

Accordingly, when factorizing $4A_{134}A_{234}$, the terms corresponding for the side $z_{3}z_{4}$ can be written as

\[ 4 A_{134}A_{234} c (a+e+d+b-2f) \]

\vspace{2 mm}

This is precisely the quantity $4 E_{34}$ as defined in \cite{M}. Repeating this computation for each of the sides of the quadrilateral, we find that:

\vspace{2 mm}

\noindent The products of sines corresponding to the side $z_{1}z_{2}$ is $4 A_{123}A_{124} f (a+b+e+d-2c)$ which is $4E_{12}$ as defined in \cite{M}.

\vspace{1 mm}

\noindent The products of sines corresponding to the side $z_{1}z_{3}$ is $-4 A_{123}A_{134} e (a+c+d+f-2b)$ which is $-4E_{13}$ as defined in \cite{M}.

\vspace{1 mm}

\noindent The products of sines corresponding to the side $z_{1}z_{4}$ is $4 A_{124}A_{134} a (c+e+b+f-2d)$ which is $4E_{14}$ as defined in \cite{M}.

\vspace{1 mm}

\noindent The products of sines corresponding to the side $z_{2}z_{3}$ is $4 A_{123}A_{234} d (c+b+e+f-2a)$ which is $4E_{23}$ as defined in \cite{M}.

\vspace{1 mm}

\noindent The products of sines corresponding to the side $z_{2}z_{4}$ is $-4 A_{124}A_{234} b (c+d+f+a-2e)$ which is $-4E_{24}$ as defined in \cite{M}.

\vspace{2 mm}

Consequently, requiring $E\geq 0$ is the same as requiring that  

\[ E_{12} + E_{23} + E_{34} + E_{14} - E_{13} - E_{24} \geq 0 \]

This is what we proved in \cite{M}. With that accomplished, we can see that the proof of the second Atiyah-Sutcliffe conjecture for a convex quadrilateral is now complete.


\begin{thebibliography}{20}
%
\bibitem{A1} M.F. Atiyah,   \emph{The geometry of classical particles}. Surveys in differential geometry $\bf{vol VII}$, 1-15, International Press, Somerville, MA  (2000)
%
\bibitem{A2}  M.F. Atiyah: \emph{Configuration of points}.  Phil. Trans. R. Lond.  $\bf{A 359}$, 1375-1387 (2001)
%

\bibitem{BR}  M.V. Berry, J.M. Robbins: \emph{Indistinguishability for quantum particles: spin, statistics and geometric phase}. Proc. Roy. Soc. London Ser. . $\bf{A 453}$,  1771-1790 (1997)
%
\bibitem{D}  D. Z. Dokovic: \emph{Verification of Atiyah's conjecture for some nonplanar configurations with dihedral symmetry}. Publ. Inst. Math. (Beograd) (N.S.) $\bf{72}$, 23-28 (2002)
%
\bibitem{EN} M. Eastwood, P. Norbury: \emph{A proof of Atiyah's conjecture on configurations of four points in Euclidean three-space}. Geom. Topol. $\bf{5}$, 885-893  (2001)
%
\bibitem{MP} M. Mazur, B. V. Penterko: \emph{On the conjectures of Atiyah and Sutcliffe}. Geom. Dedicata. DOI 10.1007/s10711-011-9636-6
%

\bibitem{AS} Michael Atiyah and Paul Sutcliffe, The geometry of point particles, R. Soc. Lond. Proc. Ser. A Math.
Phys. Eng. Sci. 458 (2002), no. 2021, 1089–1115. MR 1902577 (2003c:55019)
%

\bibitem{M1} M.  Bou Khuzam and M.  Johnson, On the conjecture regarding the 4-point Atiyah determinant,  SIGMA 10 , 070, (2014)

%
\bibitem{M} M. Bou Khuzam, "A Degree Six Inequality on Convex Quadrilaterials", http://arxiv.org/abs/2112.11765
%



\end{thebibliography}
\end{document}